\documentclass[12pt]{amsart}
\usepackage{amsbsy,amsmath,amscd,amsfonts,amsgen,amsopn,amssymb,amstext,amsthm,fullpage,setspace,indentfirst}
\usepackage[all,cmtip]{xy}
\numberwithin{equation}{section}
\newtheorem*{thma}{\textnormal{\textit{Theorem I}}}
\newtheorem*{thmb}{\textnormal{\textit{Theorem II}}}
\newtheorem*{propc}{\textnormal{\textit{Proposition III}}}
\newtheorem*{cord}{\textnormal{\textit{Corollary IV}}}
\newtheorem*{core}{\textnormal{\textit{Corollary V}}}
\newtheorem*{thmf}{\textnormal{\textit{Theorem VI}}}
\newtheorem*{thmg}{\textnormal{\textit{Theorem VII}}}
\newtheorem{thm}{\textnormal{\textit{Theorem}}}[section]
\newtheorem{prop}[thm]{\textnormal{\textit{Proposition}}}
\newtheorem{cor}[thm]{\textnormal{\textit{Corollary}}}
\newtheorem{lem}[thm]{\textnormal{\textit{Lemma}}}
\begin{document}
\title{Kernel Bundles, Syzygies Of Points, and The Effective Cone Of $C_{g-2}$}
\author{Yusuf Mustopa}
\begin{abstract}
We obtain a complete description of the effective cone of $C_{g-2}$ when $C$ is a general curve of genus $g \geq 6,$ as well as a new bound in the case where $C$ is a smooth plane quintic.  In addition, we obtain a new virtual bound for the effective cone of $C_{g-2m}$ which is a genuine bound when $m=2,$ and we also characterize certain natural divisors on $C_{g-2}$ as subordinate loci associated to adjunctions of kernel bundles. 
\end{abstract}
\maketitle
\vskip10pt
\begin{center}
\begin{tiny}
\textit{Dedicated To The Memory Of Sandra Samelson}
\end{tiny}
\end{center}
\vskip15pt
\section*{Introduction}
Let $C$ be a smooth projective curve of genus $g.$  The $d-$th symmetric power $C_{d}$ of $C$ is a smooth projective variety of dimension $d$ which is a fine moduli space parametrizing effective divisors of degree $d$ on $C;$ consequently, the ``degree-$d$ aspect" of maps from $C$ to other varieties is encoded in the geometry of $C_{d}.$   

A rich source of subvarieties of $C_{d}$ is provided by the so-called \textit{subordinate loci} associated to linear series on $C.$  Given a linear series $(L,V)$ of degree $n \geq d$ and dimension $r \leq d$ on $C,$ the subordinate locus $\Gamma_{d}(L,V)$ (we omit $V$ when the linear series is complete) is an $r-$dimensional determinantal subvariety of $C_{d}$ supported on the set $$\{D \in C_{d}: V \cap H^{0}(C,L(-D)) \neq \{0\}\}$$  These loci have been used to study the structure of algebraic cycles on the Jacobian of $C$ (e.g. \cite{Iz1}, \cite{Iz2} \cite{Herb}, \cite{vdGK}) and---more importantly for us---have also been used to study cones of divisor classes on $C_{d}$ (e.g. \cite{Kou},\cite{Pac}, \cite{Deb}, \cite{Chan}).  Our first result, which generalizes (iii) of Theorem F in \cite{Mus}, continues directly in the latter vein.        

\begin{thma}
Let $C$ be a general curve of genus $g \geq 5.$  Then the cycle $\mathfrak{D}$ supported on $$\bigcup\Gamma_{g-2}(L,V)$$ (where the union is taken over all linear series $(L,V)$ of dimension 1 and degree $g-1$) is an effective divisor which spans a boundary ray of the effective cone of $C_{g-2}.$
\end{thma}

It is not immediately clear that $\mathfrak{D}$ is of codimension 1 in $C_{g-2}.$  We will address this issue later in this introduction; for now, we point out that the linear series $(L,V)$ of dimension 1 and degree $g-1$ on a general curve $C$ of genus $g \geq 5$ are parametrized by a $(g-4)-$dimensional variety (this follows from Theorem \ref{gieseker}) and that the loci $\Gamma_{g-2}(L,V)$ are curves which cover $\mathfrak{D}$ as $(L,V)$ varies.    

Let us also mention an important property of $\mathfrak{D}$ which will be used later: whenever $D \in C_{g-2}-\mathfrak{D},$ we have that $\dim{|D|}=0$ and the residual series $|K_{C}(-D)|$ is a basepoint-free pencil (compare Corollary V). 

The notion of subordinate locus introduced above may be generalized in a straightforward manner to define subordinate loci associated to \textit{coherent systems} on $C,$ that is, pairs $(F,V)$ which consist of a vector bundle $F$ on $C$ and a subspace $V \subseteq H^{0}(C,F).$  One of the main goals of this paper is to demonstrate how these higher-rank subordinate loci occur in the study of the effective cone of $C_{d}.$      

Our results on these loci center around the well-known \textit{kernel bundle} $M_{L}:=\phi_{L}^{\ast}\Omega^{1}_{\mathbb{P}^{N}}(1),$ where $L$ is a globally generated line bundle on $C$ and $\phi_{L}:C \rightarrow \mathbb{P}^{N}$ is the associated morphism.   

\begin{thmb}
Let $C$ be a general curve of genus $g \geq 5.$  Then for all $p_{0} \in C$ we have the equality of cycles 
\begin{equation}
\label{cyceq}
\Gamma_{g-2}(K_{C} \otimes M_{K_{C}(-p_{0})})=\mathfrak{D}+X_{p_{0}}
\end{equation}
where $X_{p_{0}}$ is the reduced divisor supported on $\{D'+p_{0}:D' \in C_{g-3}\}.$
\end{thmb}  
\noindent

Subordinate loci of the form $\Gamma_{d}(K_{C} \otimes M_{L})$ admit a natural interpretation in terms of multiplication maps.  Pulling back the Euler sequence on $\mathbb{P}^{N}$ via $\phi_{L}$ yields the exact sequence
\begin{equation}
\label{euler}
0 \rightarrow M_{L} \rightarrow H^{0}(C,L) \otimes \mathcal{O}_{C} \rightarrow L \rightarrow 0
\end{equation}
Fix an element $D \in C_{d}.$  Twisting by $K_{C}(-D)$ and taking cohomology yields the exact sequence 
$$
\begin{tiny}
\begin{CD}
0 @>>> H^{0}(C,K_{C} \otimes M_{L}(-D)) @>>> H^{0}(C,L) \otimes H^{0}(C,K_{C}(-D)) @>{\mu_{L,K_{C}(-D)}}>> H^{0}(C,K_{C} \otimes L(-D))\\
\end{CD}
\end{tiny}$$
where $\mu_{L,K_{C}(-D)}$ is the multiplication map.  It then follows that $\Gamma_{d}(K_{C} \otimes M_{L})$ parametrizes all $D \in C_{d}$ for which $\mu_{L,K_{C}(-D)}$ fails to be injective.  

Any serious consideration of multiplication maps is closely related to Koszul cohomology, and the reader familiar with the latter has undoubtedly noticed its relevance to the present work.  While we make almost no explicit use of it here, it is almost certainly the natural framework for further progress along some of the lines indicated by our results on higher-rank subordinate loci.  We refer to \cite{Gr} for a foundational treatment of Koszul cohomology and to \cite{ApFa} for a nice survey of some recent developments.       

Let us now return to the issues surrounding Theorem I.  We begin by introducing the basic divisor classes on $C_{d}.$  For each $p \in C$ the image $X_{p}$ of the embedding $$i_{p}:C_{d-1} \hookrightarrow C_{d}, \hspace{0.1cm} D' \mapsto D'+p$$ is an ample divisor on $C_{d}$ whose numerical class (which is independent of $p$) will be denoted by $x.$  The pullback of a theta-divisor on $\textrm{Pic}^{d}(C)$ via the natural map $$a_{d}:C_{d} \rightarrow \textrm{Pic}^{d}(C), \hspace{0.1cm} D \mapsto \mathcal{O}_{C}(D)$$ is a nef divisor on $C_{d}$ whose numerical class will be denoted by $\theta.$  (When $2 \leq d \leq g,$ we also have that $\theta$ is big.)  The classes $x$ and $\theta$ are linearly independent in the real N\'{e}ron-Severi space $N^{1}_{\mathbb{R}}(C_{d}),$ and when $C$ is general $N^{1}_{\mathbb{R}}(C_{d})$ is generated by $x$ and $\theta.$     

The following result plays an important role in the proofs of Theorems I and II for reasons that will soon be evident.  It also gives a new bound for the effective cone of $C_{g-4}.$

\begin{propc}
\label{classes}
Let $C$ be a general curve of genus $g \geq 5,$ and let $m$ be a positive integer which is at most $\frac{g}{2}-1.$  Then there is a virtual divisor $\mathfrak{D}_{m}$ on $C_{g-2m}$ supported on the set $$\bigcup\Gamma_{g-2m}(\mathcal{L},V)$$ (where the union is taken over all linear series $(\mathcal{L},V)$ of degree $g-m$ and dimension 1) with class proportional to $\theta-(1+\frac{2m}{g-2m})x,$ and it is an effective divisor when $m=1$ or $m=2.$  
\end{propc} 

We now proceed to list some consequences in the case $m=1.$  Theorem 3 of \cite{Kou} says that for all $C$ and for all $d \geq 2,$ the class $$\Delta=2(-\theta+(g+d-1)x)$$ of the diagonal locus on $C_{d}$ (which parametrizes effective divisors of degree $d$ having multiplicity) spans a boundary ray of the effective cone of $C_{d}.$  Combining this result with Theorem I and Proposition III yields 

\begin{cord}
\label{effchar}
If $C$ is a general curve of genus $g \geq 5$ then the effective cone of $C_{g-2}$ is spanned by the classes $\Delta$ and $\theta-(1+\frac{2}{g-2})x.$ 
\end{cord}

Recall that for nonnegative integers $s$ and $e$ there exists a fine moduli variety $G^{s}_{e}(C)$ parametrizing linear series of degree $e$ and dimension $s$ on $C.$  In \cite{Mus} we made heavy use of the fact that for a general curve $C$ of genus $g \geq 3$ and $2 \leq d \leq g-1,$ the rational map $\tau:C_{d} \dashrightarrow G^{g-d-1}_{2g-2-d}(C)$ taking a general $D \in C_{d}$ to the residual linear series $|K_{C}(-D)|$ is an isomorphism of smooth varieties in codimension 1.  In particular, $\tau$ preserves the effective cone, so that Corollary IV implies        

\begin{core}
If $C$ is a general curve of genus $g \geq 5,$ then the effective cone of $G^{1}_{g}(C)$ is spanned by class of the divisor $\mathfrak{D}'$ supported on $$\{(L,V) \in G^{1}_{g}(C): (L,V) \textnormal{ has a base point}\}$$ and the class of the divisor ${\Delta}'$ supported on the closure of $$\{(L,V) \in G^{1}_{g}(C): V=H^{0}(L)\textnormal{ and }|K_{C} \otimes L^{-1}|\textnormal{ has multiplicity}\}$$  
\end{core}

Let us now consider the case $m \geq 2,$ where $\mathfrak{D}_{m}$ admits a nice interpretation in terms of secant planes.  It can be shown that for a general curve $C$ and general $D \in C_{g-2m},$ the residual series $|K_{C}(-D)|$ is basepoint-free of dimension $2m-1$ and the associated morphism $\phi_{K_{C}(-D)}$ is birational onto its image (we refer to the proof of Proposition III for details).  The Riemann-Roch Theorem then implies that for $m \geq 2,$ the support of $\mathfrak{D}_{m}$ may be characterized as the Zariski closure of the set $$\{D \in C_{g-2m}:\dim{|K_{C}(-D)|}=2m-1\textnormal{ and }\phi_{K_{C}(-D)}(C)\textnormal{ admits an }m-\textnormal{secant }(m-2)-\textnormal{plane}\}$$  The structure of $\mathfrak{D}_{m}$ as a virtual divisor on $C_{g-2m}$ comes from the fact that it is supported on the projection of a locus in $C_{g-2m} \times C_{m}$ which is shown to be of pure dimension $g-2m-1,$ namely $$\{(D,E) \in C_{g-2m} \times C_{m}: E\textnormal{ fails to impose independent conditions on }|K_{C}(-D)|\}$$  Showing that $\mathfrak{D}_{m}$ is an honest divisor is therefore equivalent to showing that for general $D \in C_{g-2m},$ there are at most finitely many $E \in C_{m}$ which fail to impose independent conditions on $|K_{C}(-D)|.$   For recent progress on problems of this nature, see \cite{Cot}, \cite{Fa}, and the references therein.

We now turn from general curves to special curves.  In the hyperelliptic case, the non-diagonal boundary ray of the effective cone of $C_{g-2}$ (resp. $C_{g-1}$) in the $(\theta,x)-$plane is spanned by $\theta-3x$ (resp. $\theta-2x$); this is part of Proposition H in \cite{Mus}.  In the trigonal case it is spanned by $\theta-2x$ (this is easily extracted from the proof of Theorem 5 in \cite{Kou}).  These special families, together with smooth plane quintics, occupy important places in two distinct classification schemes--namely the Martens-Mumford classification of Brill-Noether loci and classification of curves by their Clifford index (which we will describe shortly). 

Our proof of Theorem I cannot be extended to the smooth plane quintic case.  This is essentially due to the fact that every complete linear series of degree 5 and dimension 1 on a smooth plane quintic has a base point.  However, we have the following pleasant characterization of $\mathfrak{D}:$

\begin{thmf}
Let $C$ be a smooth plane quintic, and let $\mathcal{Q}=\mathcal{O}_{C}(1).$  Then we have the equality of cycles
\begin{equation}
\label{quin}
\mathfrak{D}=2 \cdot \Gamma_{4}(K_{C} \otimes M_{\mathcal{Q}}).   
\end{equation}
\end{thmf}  
The overt syzygetic content in this paper belongs to the proof of this theorem, which employs the minimal graded resolution of 4 points in $\mathbb{P}^{2}$ in a crucial way.   

We also obtain a strict outer bound for the effective cone of $C_{4}.$

\begin{thmg}
If $C$ is a smooth plane quintic, then there is no effective divisor on $C_{4}$ whose class is proportional to $\theta-2x.$
\end{thmg} 

We conclude this introduction with a brief discussion of the Clifford index of a curve and its (possible) relation to our work.  We define $\textnormal{Cliff}(C)$ to be $$\min\{\deg{L}-2\dim{|L|}: L \textnormal{ a line bundle on }C \ni \dim{|L|}\geq 1, \dim{|K_{C} \otimes L^{-1}|} \geq 1\}.$$  This is a rough measure of how general $C$ is in the sense of moduli; one can check, for instance, that $\textnormal{Cliff}(C)=0$ precisely when $C$ is hyperelliptic and that $\textnormal{Cliff}(C)=1$ precisely when $C$ is either trigonal or a smooth plane quintic.  We mention in passing that the bundle $M_{K_{C}}$ has been intensively studied with an eye towards the Green conjecture relating $\textnormal{Cliff}(C)$ to the syzygies of the canonical embedding of $C$ (e.g. \cite{PR},\cite{Ein}) and that the tools developed by Voisin \cite{V} in her proof of this conjecture for a general curve of even genus $g$ have been used by Pacienza \cite{Pac} to compute the nef cone of $C_{\frac{g}{2}}$ for all such curves $C.$         

Given that the effective cones of $C_{g-1}$ and $C_{g-2}$ detect whether $\textnormal{Cliff}(C)$ is 0 or greater, it is natural in light of the trigonal case and Theorem VII to ask if the effective cone of $C_{g-2}$ detects whether $\textnormal{Cliff}(C)$ is 1 or greater.  Theorem VII does not preclude the possibility that the effective cone of $C_{4}$ has an open boundary ray spanned by $\theta-2x$ when $C$ is a smooth plane quintic, and this points to two extreme scenarios: one where Theorem I continues to hold in the plane quintic case, and another where the closure of the part of the effective cone of $C_{g-2}$ in the $(\theta,x)-$plane is the same for all curves $C$ satisfying $\textnormal{Cliff}(C)=1.$  As of this writing, however, we have yet to find evidence that either of these is true.               

\vskip15pt

\textit{Notation and Conventions:}  We work over the field of complex numbers.  $C$ will always denote a smooth projective curve.  The symbol $\mathfrak{g}^{r}_{d}$ is used to refer to a linear series (complete or otherwise) on $C$ of degree $r$ and dimension $d.$  All cycle classes on smooth varieties lie in the algebraic cohomology ring with coefficients in $\mathbb{R}$.  We say that a property holds for \textit{a general curve of genus $g$} if it holds on the complement of the union of countably many proper subvarieties of the moduli space $\mathcal{M}_{g}$ of smooth projective curves of genus $g$.

\section{Preliminaries on $C_{d}$}
\subsection{\textnormal{\textit{Intersection Theory}}}
The following formula, which will be used freely, is a consequence of the Poincar\'{e} formula (p.25 of \cite{ACGH}).

\begin{lem}
\label{intersect}
For all $0 \leq k \leq d \leq g,$ $$x^{k}\theta^{d-k}=\displaystyle\frac{g!}{(g-d+k)!}$$
\end{lem}

\subsection{\textnormal{\textit{Class Formulas}}}
We begin with the basic result on classes of subordinate loci.  This is Lemma 3.2 on p.342 of \textit{loc. cit.}

\begin{lem}
\label{classcomp}
Let $(L,V)$ be a $\mathfrak{g}^{r}_{n}$ on $C,$ and assume $r \leq d \leq n.$  Then $\Gamma_{d}(L,V)$ is an $r-$dimensional subvariety of $C_{d}$ whose class is $$\gamma_{d}(\mathfrak{g}^{r}_{n}):=\displaystyle\sum_{j=0}^{d-r}\displaystyle\binom{n-g-r}{j}\displaystyle\frac{x^{j} \cdot \theta^{d-r-j}}{(d-r-j)!}$$
\end{lem}

While the following result is a corollary in the sense that it can be obtained in a straightforward manner from Lemma \ref{classcomp}, it can be obtained much more easily from the basic properties of the morphism $a_{d}.$    

\begin{cor}
\label{exc}
$\theta \cdot \gamma_{d}(\mathfrak{g}^{1}_{d})=0$ and $x \cdot \gamma_{d}(\mathfrak{g}^{1}_{d})=1.$ \hfill \qedsymbol 
\end{cor}
The next two formulas are crucial to the proofs of Proposition III and Theorem I.  Recall that $C^{1}_{d},$ which is the exceptional locus of the map $a_{d},$ parametrizes effective divisors of degree $d$ that move in a positive-dimensional linear series.  The next formula follows immediately from the Theorem on p.326 of \textit{loc. cit.}
 
\begin{prop}
\label{seeonedee}
If the cycle $C^{1}_{d}$ is either empty or of the expected dimension $2d-g-1,$ its class $c^{1}_{d}$ is equal to $$\displaystyle\frac{\theta^{g-d+1}}{(g-d+1)!}-\displaystyle\frac{x\theta^{g-d}}{(g-d)!}$$   
\end{prop}

\begin{prop}
\label{pushpull}
(Push-Pull formula) For each algebraic cycle $Z$ in $C_{d},$ let $B_{k}(Z)$ denote the class of the cycle parametrizing the set $\{E \in C_{d-k} : |D-E| \neq \emptyset \textnormal{ for some }D \in Z\}.$  Then we have the formula $$B_{k}(x^{\alpha}\theta^{\beta})=\displaystyle\sum_{j=0}^{k}\displaystyle\binom{\alpha}{k-j}\displaystyle\binom{\beta}{j}\displaystyle\binom{g-\beta+j}{j}j!x^{\alpha-k+j}\theta^{\beta-j} $$
\end{prop}

The derivation of the Push-Pull formula is outlined in exercise batch D of Chapter VIII of \textit{loc. cit.}

\subsection{\textnormal{\textit{Special Subordinate Loci Associated To Basepoint-Free Pencils}}}
The classical subordinate loci of primary interest to us are of the form $\Gamma_{d}(L,V)$ where $(L,V)$ is a $\mathfrak{g}^{1}_{d+1}.$  

\begin{lem}
\label{genus}
Let $C$ be a curve of genus $g \geq 1,$ and let $(L,V)$ be a basepoint-free $\mathfrak{g}^{1}_{d+1}$ on $C$.  Then $\Gamma_{d}(L,V) \cong C.$
\end{lem}

\begin{proof}
(i) Let $\mathfrak{Y}$ be the incidence variety supported on $\{(p,D) \in C \times C_{d}:p+D \in |V|\}.$  We will proceed by showing that $C$ and $\Gamma_{d}(L,V)$ are both isomorphic to $\mathfrak{Y}.$      

Since $(L,V)$ is basepoint-free, we have well-defined bijective morphisms $j_{1}:C \rightarrow \mathfrak{Y}$ and $j_{2}:\Gamma_{d}(L,V) \rightarrow \mathfrak{Y}$ given by $j_{1}(p)=(p,|V| \cap |L(-p)|)$ and $j_{2}(D)=(|L(-D)|,D).$  This implies that $g \leq p_{a}(\mathfrak{Y})$ and $p_{a}(\Gamma(L,V)) \leq p_{a}(\mathfrak{Y}).$

The morphisms $j_{1}$ and $j_{2}$ are sections of the projections $\rho_{1}:\mathfrak{Y} \rightarrow C$ and $\rho_{2}:\mathfrak{Y} \rightarrow \Gamma_{d}(L,V),$ respectively.  Consequently $\rho_{1}$ and $\rho_{2}$ are bijective, so that $p_{a}(\mathfrak{Y}) \leq g$ and $p_{a}(\mathfrak{Y}) \leq p_{a}(\Gamma_{d}(L,V)).$  It follows that the arithmetic genera of the curves $C, \Gamma_{d}(L,V),$ and $\mathfrak{Y}$ are all equal; this implies that $j_{1}$ and $j_{2}$ are isomorphisms.     
\end{proof}

\subsection{\textnormal{\textit{Higher-Rank Subordinate Loci}}}
\label{higher}
Let $d \geq 2$ be given, let $\pi_{1},\pi_{2}$ be the projection maps associated to $C \times C_{d},$ and let $\mathcal{U}$ be the universal divisor on $C \times C_{d},$ i.e. the incidence variety supported on $\{(p,D) \in C \times C_{d}: p \in \textnormal{supp}(D)\}.$  For any vector bundle $F$ of rank $r$ on $C,$ we define $$\mathbb{E}(F):=\pi_{2\ast}(\pi^{\ast}_{1}F \otimes \mathcal{O}_{\mathcal{U}})$$ which is a vector bundle of rank $rd$ on $C_{d}$ whose fibre over $D \in C_{d}$ is $H^{0}(D,F|_{D}).$  If $(F,V)$ is a coherent system, the restriction maps $V \rightarrow H^{0}(D,F|_{D})$ globalize to a vector bundle morphism $$\alpha_{(F,V)}:V \otimes \mathcal{O}_{C_{d}} \rightarrow \mathbb{E}(F).$$  

\noindent
\textit{Definition:}  $\Gamma_{d}(F,V)$ is the degeneracy locus associated to $\alpha_{(F,V)}.$  If $V=H^{0}(C,F),$ we will denote $\Gamma_{d}(F,V)$ by $\Gamma_{d}(F).$
 
\begin{prop}
\label{chern}
Let $(F,V)$ be a coherent system on $C$ of rank $r \geq 2$ and degree $f$ such that $rd=\dim{V}.$  Then the class of the virtual divisor $\Gamma_{d}(F,V)$ on $C_{d}$ is $r\theta-(rd+rg-f-r)x.$
\end{prop}
\begin{proof}
By Porteous' Formula, the class we seek is $c_{1}(\mathbb{E}(F)),$ which we will extract from the Chern character of $\mathbb{E}(F).$  

The calculation is almost identical to the proof of Lemma 2.5 on p.340 of \cite{ACGH}.  In what follows, $\eta$ denotes the class of the pullback via $\pi_{1}$ of a point on $C$ and the pullbacks via $\pi_{2}$ of the classes $x$ and $\theta$ are also denoted by $x$ and $\theta,$ respectively.  By Grothendieck-Riemann-Roch, we have that $$\textnormal{td}(C_{d}) \cdot \textnormal{ch}(\mathbb{E}(F))=(\pi_{2})_{\ast}(\textnormal{td}(C \times C_{d}) \cdot \textnormal{ch}(\pi_{1}^{\ast}F \otimes \mathcal{O}_{\mathcal{U}}))$$ Canceling out the $\textnormal{td}(C_{d})$ term and using the short exact sequence $$0 \rightarrow \pi_{1}^{\ast}F(-\mathcal{U}) \rightarrow \pi_{1}^{\ast}F \rightarrow \pi_{1}^{\ast}F \otimes \mathcal{O}_{\mathcal{U}} \rightarrow 0$$ together with the formula for the class of $\mathcal{U}$ from p.338 of \textit{loc. cit.}, we have that $$\textnormal{ch}(\mathbb{E}(F))=(\pi_{2})_{\ast}((1+(1-g)\eta) \cdot ((r+f\eta)-(r+(f-rd)\eta-r\eta\theta)e^{-x}))$$ $$=(\pi_{2})_{\ast}(r+(f+r(1-g))\eta-(r+(f-rd+r-rg)\eta-r\eta\theta)e^{-x})$$ $$=(f+r(1-g))+(rd+rg-f-r+r\theta)e^{-x}.$$
\end{proof}

\section{Proofs of Proposition III and Theorem I}
\label{thema}
We will require two major results of Brill-Noether theory.  The first of these is Gieseker's Theorem (e.g. Theorem 1.6 on p. 214 of \cite{ACGH}):

\begin{thm} 
\label{gieseker}
Let $C$ be a general curve of genus $g.$  Let $d,r$ be integers satisfying $d \geq 1$ and $r \geq 0.$  Then $G^{r}_{d}(C)$ is smooth of dimension $g-(r+1)(g-d+r).$ 
\end{thm}

We will also use the following theorem of Fulton and Lazarsfeld (e.g. Theorem 1.4 on p.212 of \textit{loc. cit.}): 

\begin{thm}
\label{fullaz}
Let $C$ be a smooth curve of genus $g.$  Let $d,r$ be integers such that $d \geq 1, r \geq 0.$ Assume that $$g-(r+1)(g-d+r) \geq 1.$$  Then $G^{r}_{d}(C)$ is connected.
\end{thm}

The following corollaries are immediate applications of these theorems; the first follows from Theorem \ref{gieseker}, while the second follows from combining Theorems \ref{gieseker} and \ref{fullaz}. 

\begin{cor}
\label{base}
If $C$ is a general curve of genus $g \geq 4,$ then the basepoint-free $\mathfrak{g}^{1}_{g-1}$'s are a nonempty Zariski open subset of $G^{1}_{g-1}(C).$ \hfill \qedsymbol
\end{cor}

\begin{cor}
\label{gifl}
If $C$ is a general curve of genus $g,$ and $d,r$ are integers such that $d \geq 1, r \geq 0,$ and $g-(r+1)(g-d+r) \geq 1,$ then $G^{r}_{d}(C)$ is irreducible. \hfill \qedsymbol
\end{cor}

\textit{Proof of Proposition III:}  Let $C$ be a general curve of genus $g \geq 5,$ and let $m$ be a positive integer satisfying $m \leq \frac{g}{2}-1.$  Consider the diagram  

$$\xymatrix{
C_{g-2m} \times C_{m} \ar[d]^{\pi} \ar[r]^{\hspace{0.6cm} \sigma} &C_{g-m}\\
C_{g-2m}}$$ where $\sigma$ is the addition map and $\pi$ is projection onto the first factor.  Since $C$ is general, the Brill-Noether Theorem implies that $C^{1}_{g-m}$ is of pure dimension $g-2m-1,$ and since $\sigma$ is a finite map of smooth varieties, the same is true of $\sigma^{-1}(C^{1}_{g-m}).$  Therefore $\mathfrak{D}_{m}:=\pi_{\ast}\sigma^{\ast}(C^{1}_{g-m})$ is a virtual divisor on $C_{g-2m}$.

We now compute the class of $\mathfrak{D}_{m}.$  We have from Proposition \ref{seeonedee} that $$c^{1}_{g-m}=\displaystyle\frac{\theta^{m+1}}{(m+1)!}-\displaystyle\frac{x\theta^{m}}{m!}$$  Combining this with Proposition \ref{pushpull}, we have that the class of $\mathfrak{D}_{m}$ is  

\begin{equation}
\label{impclass}
B_{m}(c^{1}_{g-m})=\displaystyle\binom{g}{m}\bigg(\displaystyle\frac{g-2m}{g} \cdot \theta-x\bigg).
\end{equation}

The cycle $\mathfrak{D}_{m}$ is an effective divisor on $C_{g-2m}$ precisely when the intersection $\sigma^{-1}(C^{1}_{g-m}) \cap \pi^{-1}(D)$ is at most finite for general $D \in C_{g-2m}.$  Our next and final task is to verify this condition for $m=1$ and $m=2.$  Note that since $\sigma^{-1}(C^{1}_{g-m}) \cap \pi^{-1}(D)=\{D\} \times C_{m}$ for all $D \in C^{1}_{g-2m},$ we must necessarily consider a general element of $C_{g-2m}-C^{1}_{g-2m}.$  (The Brill-Noether Theorem implies that $C^{1}_{g-2m}$ is a proper Zariski-closed subset of $C_{g-2m}$ for all $m \geq 1.$)  

If $m=1,$ then $\sigma^{-1}(C^{1}_{g-1}) \cap \pi^{-1}(D)$ parametrizes the points $q \in C$ for which $\dim{|D+q|} \geq 1,$ or equivalently, the base points of $|K_{C}(-D)|.$  If $D \in C_{g-2}-C^{1}_{g-2},$ then $|K_{C}(-D)|$ is 1-dimensional, and thus has at most finitely many base points.  

On our way to settling the case $m=2,$ we justify the claim (stated in the Introduction) that for all $m \geq 2$ and general $D \in C_{g-2m},$ the linear series $|K_{C}(-D)|$ is basepoint-free of dimension $2m-1$ and the associated morphism $\phi_{K_{C}(-D)}$ is birational onto its image.  We assume throughout that $D \in C_{g-2m}-C^{1}_{g-2m}.$  By Riemann-Roch, we have that $\dim{|K_{C}(-D)|}=2m-1,$ and by Brill-Noether, the set of all $D \in C_{g-2m}$ for which $|K_{C}(-D)|$ admits a base point is a proper Zariski-closed subset of $C_{g-2m};$ this establishes the first part of the claim.  If $\phi_{K_{C}(-D)}$ is not birational onto its image, its image is necessarily a nondegenerate rational curve of degree at least $2m-1$ in $\mathbb{P}^{2m-1}$ since a general curve of genus $g$ cannot be a multiple cover of an irrational curve (cf. Exercise C-6 in Chapter VIII of \cite{ACGH}).  The morphism $\phi_{K_{C}(-D)}$ then factors through a branched cover of $\mathbb{P}^{1}$ by $C$ having degree at most $\frac{g+2m-2}{2m-1}.$  But this is impossible, since any branched cover of $\mathbb{P}^{1}$ by $C$ must have degree at least $\frac{g}{2}+1$ by Brill-Noether.      

The characterization of $\mathfrak{D}_{m}$ in terms of secant planes stated in the Introduction implies that $\mathfrak{D}_{2}$ parametrizes $D \in C_{g-4}$ for which $\phi_{K_{C}(-D)}$ fails to be an embedding.  Since $\phi_{K_{C}(-D)}$ is birational onto its image for general $D \in C_{g-4}$ and the image of $\phi_{K_{C}(-D)}$ has at most finitely many singular points, it follows that $\sigma^{-1}(C^{1}_{g-2}) \cap \pi^{-1}(D)$ is at most finite for general $D \in C_{g-4}.$  This concludes the proof that $\mathfrak{D}_{2}$ is an effective divisor. \hfill \qedsymbol

\medskip

\textit{Proof of Theorem I:}  We begin by showing the irreducibility of $\mathfrak{D},$ which we know to be of codimension 1 in $C_{g-2}$ by Proposition III.  There is a natural rational map $$\phi: C \times G^{1}_{g-1}(C) \dashrightarrow \mathfrak{D}, \enspace (p,\mathcal{L}) \mapsto |\mathcal{L}(-p)|$$ which is well-defined (outside a closed set) and dominant thanks to Corollary \ref{base}.  By Corollary \ref{gifl}, $C \times G^{1}_{g-1}(C)$ is irreducible, so $\mathfrak{D}$ is irreducible as well.

Let us now consider the curves $\Gamma_{g-2}(L,V),$ where $(L,V)$ is a $\mathfrak{g}^{1}_{g-1}$ on $C.$  By Lemma \ref{classcomp}, we have that their common numerical class is  $$\gamma_{g-2}(\mathfrak{g}^{1}_{g-1})=\displaystyle\sum_{j=0}^{g-3}\displaystyle\binom{-2}{j}\displaystyle\frac{x^{j}\theta^{g-3-j}}{(g-3-j)!}$$  It then follows that 
$$\theta \cdot \gamma_{g-2}(\mathfrak{g}^{1}_{g-1})=\displaystyle\sum_{j=0}^{g-3}(-1)^{j}(j+1)\displaystyle\frac{g!}{(j+2)!(g-3-j)!}$$ $$=g \cdot \displaystyle\sum_{j=0}^{g-3}(-1)^{j}(j+1)\displaystyle\binom{g-1}{j+2}=g \cdot \bigg(\displaystyle\sum_{j=0}^{g-3}(-1)^{j}(j+2)\displaystyle\binom{g-1}{j+2}-\displaystyle\sum_{j=0}^{g-3}(-1)^{j}\displaystyle\binom{g-1}{j+2}\bigg)$$ $$=g \cdot \bigg(\displaystyle\sum_{j=2}^{g-1}(-1)^{j}j\displaystyle\binom{g-1}{j}-\displaystyle\sum_{j=2}^{g-1}(-1)^{j}\displaystyle\binom{g-1}{j}\bigg)=g \cdot ((g-1)-(g-2))=g.$$  We need one other intersection number: $$x \cdot \gamma_{g-2}(\mathfrak{g}^{1}_{g-1})=\displaystyle\sum_{j=0}^{g-3}(-1)^{j}(j+1)\displaystyle\frac{g!}{(j+3)!(g-3-j)!}=\displaystyle\sum_{j=0}^{g-3}(-1)^{j}(j+1)\displaystyle\binom{g}{j+3}$$ $$=-\displaystyle\sum_{j=3}^{g}(-1)^{j}(j-2)\displaystyle\binom{g}{j}=2\displaystyle\sum_{j=3}^{g}(-1)^{j}\displaystyle\binom{g}{j}-\displaystyle\sum_{j=3}^{g}(-1)^{j}j\displaystyle\binom{g}{j}=g-2.$$  It then follows from setting $m=1$ in (\ref{impclass}) that $$\mathfrak{D} \cdot \gamma_{g-2}(\mathfrak{g}^{1}_{g-1})=((g-2)\theta-gx) \cdot \gamma_{g-2}(\mathfrak{g}^{1}_{g-1})=0.$$  Suppose there exists an irreducible effective divisor $\mathfrak{E}$ on $C_{g-2}$ whose class is proportional to $\theta-tx$ for some $t > 1+\frac{2}{g-2}.$  Then $\mathfrak{E} \cdot \gamma_{g-2}(\mathfrak{g}^{1}_{g-1}) < 0.$  Since $\Gamma_{g-2}(L,V)$ is irreducible for general $(L,V)$ by Lemma \ref{genus} and Corollary \ref{base}, we have that the curves $\Gamma_{g-2}(L,V)$ cover $\mathfrak{E}.$  We may then conclude that $\mathfrak{D}=\mathfrak{E},$ which is absurd. \hfill \qedsymbol      

\section{The Virtual Divisor $\Gamma_{d}(K_{C} \otimes M_{L})$}
In preparation for the proofs of Theorems II and VI, we present two results that yield two different ways of showing that $\Gamma_{d}(K_{C} \otimes M_{L})$ is a virtual divisor in the cases of interest to us.  

The first result, which we state without proof, is a straightforward consequence of applying Green's ``$K_{p,1}-$Theorem" (Theorem 3.c.1 in \cite{Gr}) and Lazarsfeld's description of Koszul cohomology in terms of kernel bundles (e.g. Theorem 2.6 in \cite{ApFa}) to the Koszul cohomology group $K_{r-1,r}(C,L)$.  It can also be viewed as a slight modification of Proposition 3.4 in \cite{Tei}.  

\begin{prop}
\label{bpf}
Let $C$ be an irrational curve and let $k \geq 4$ be an integer such that the general complete pencil of degree $k$ on $C$ is basepoint-free.  Let $L$ be a line bundle of degree $2k-3$ such that $h^{0}(C,L)=k-1.$  Then $h^{0}(C,M^{\ast}_{L})=k-1.$ \hfill \qedsymbol 
\end{prop}

For the statement and proof of the second result, we reuse the notation from Section \ref{higher}.

\begin{prop}
\label{mult}
Let $C$ be a nonhyperelliptic curve, let $d$ be an integer satisfying $2 \leq d \leq g-1,$ and define $C^{\ast}_{d}$ to be the quasiprojective variety $C_{d}-C^{1}_{d}.$  Furthermore, define $\mathcal{U}^{\ast}$ to be the restriction of $\mathcal{U}$ to $C \times C^{\ast}_{d},$ let $r$ be an integer satisfying $r \geq \frac{d}{g-d},$ and let $L$ be a line bundle on $C$ of degree $r(g-d)+1.$  Then the degeneracy locus of the multiplication map $$\mu_{L}:H^{0}(C,L) \otimes \pi_{1\ast}(\pi_{2}^{\ast}K_{C} \otimes \mathcal{O}(-\mathcal{U}^{\ast})) \rightarrow \pi_{1\ast}(\pi_{2}^{\ast}(K_{C} \otimes L) \otimes \mathcal{O}(-\mathcal{U}^{\ast}))$$ on $C^{\ast}_{d}$ extends to a virtual divisor on $C_{d}$ whose class is $r\theta-(r+1)x.$   
\end{prop}

\begin{proof}
For simplicity we define $\mathfrak{F}:=\pi_{1\ast}(\pi_{2}^{\ast}K_{C} \otimes \mathcal{O}(-\mathcal{U}^{\ast}))$ and $\mathfrak{G}_{L}:=\pi_{1\ast}(\pi_{2}^{\ast}(K_{C} \otimes L) \otimes \mathcal{O}(-\mathcal{U}^{\ast})).$  By Riemann-Roch and Grauert's Theorem, $\mathfrak{F}$ and $\mathfrak{G}_{L}$ are locally free sheaves on $C^{\ast}_{d}$ of respective ranks $g-d$ and $(r+1)(g-d).$  As $D$ varies over $C^{\ast}_{d},$ the exact sequence $$0 \rightarrow H^{0}(C,K_{C}(-D)) \rightarrow H^{0}(C,K_{C}) \rightarrow H^{0}(D,K_{C}|_{D}) \rightarrow 0$$ determines the following exact sequence of locally free sheaves on $C^{\ast}_{d}:$ $$0 \rightarrow \mathfrak{F} \rightarrow H^{0}(C,K_{C}) \otimes \mathcal{O}_{C^{\ast}_{d}} \rightarrow \Omega^{1}_{C_{d}} \otimes \mathcal{O}_{C^{\ast}_{d}} \rightarrow 0$$  It follows from taking determinants that $\det\mathfrak{F}=K^{-1}_{C_{d}} \otimes \mathcal{O}_{C^{\ast}_{d}}.$  Since $C$ is nonhyperelliptic, we have by Martens' Theorem that $C^{1}_{d}$ is of codimension 2 or greater.  Consequently $\det\mathfrak{F}$ extends over $C^{1}_{d}$ to $K^{-1}_{C_{d}}$ by Hartogs' Theorem, and $$c_{1}(K^{-1}_{C_{d}})=-\theta-(g-d-1)x.$$  A similar argument shows that $\det{\mathfrak{G}_{L}}$ extends over $C^{1}_{d}$ to $\det{\mathbb{E}(K_{C} \otimes L)^{\ast}},$ so that $$c_{1}(\mathfrak{G}_{L})=-c_{1}(\Gamma_{d}(K_{C} \otimes L))=-\theta-(r+1)(g-d)x$$ by Lemma \ref{classcomp} or Proposition \ref{chern}.  Therefore the class of our degeneracy locus is $$c_{1}(\mathfrak{G}_{L})-(r+1)c_{1}(K^{-1}_{C_{d}})=r\theta-(r+1)x.$$ 
\end{proof}

\section{Proof of Theorem II}

The following is an immediate consequence of either Proposition \ref{mult} or a combination of Propositions \ref{chern} and \ref{bpf}.

\begin{lem}
\label{minuspt}
If $C$ is nonhyperelliptic of genus $g \geq 5$, then for all $p_{0} \in C$ we have that $\Gamma_{g-2}(K_{C} \otimes M_{K_{C}(-p_{0})})$ is a virtual divisor on $C_{g-2}$ with class $(g-2)\theta-(g-1)x.$ \hfill \qedsymbol
\end{lem}  

\textit{Proof of Theorem II:} We begin by simultaneously showing that the virtual divisor $\Gamma_{g-2}(K_{C} \otimes M_{K_{C}(-p_{0})})$ is an honest divisor and that it contains $X_{p_{0}}.$  Let $D \in C_{g-2}-\mathfrak{D}$ be given.  Then the degree-$g$ pencil $|K_{C}(-D)|$ is basepoint-free, and we have the short exact sequence $$0 \rightarrow T_{C}(D) \rightarrow H^{0}(C,K_{C}(-D)) \otimes \mathcal{O}_{C} \rightarrow K_{C}(-D) \rightarrow 0.$$  Twisting by $K_{C}(-p_{0})$ and taking cohomology, we see that the kernel of $\mu_{K_{C}(-p_{0}),K_{C}(-D)}$ is isomorphic to $H^{0}(C,\mathcal{O}_{C}(D-p_{0})).$  We then have a natural isomorphism $$H^{0}(C,K_{C} \otimes M_{K_{C}(-p_{0})}(-D)) \cong H^{0}(C,\mathcal{O}_{C}(D-p_{0})).$$  Recall from the Introduction that $\Gamma_{g-2}(K_{C} \otimes M_{K_{C}(-p_{0})})$ is supported on the set of all $D \in C_{g-2}$ for which $\mu_{K_{C}(-p_{0}),K_{C}(-D)}$ fails to be injective.  It then follows that $D \in C_{g-2}-\mathfrak{D}$ is contained in $\Gamma_{g-2}(K_{C} \otimes E^{\ast}_{K_{C}(-p_{0})})$ precisely when $D \in X_{p_{0}}.$  

Since $\mathfrak{D}$ is irreducible and $\Gamma_{g-2}(K_{C} \otimes M_{K_{C}(-p_{0})})$ is a determinantal subscheme of $C_{g-2}$ having the expected dimension, either $\Gamma_{g-2}(K_{C} \otimes M_{K_{C}(-p_{0})})$ contains all of $\mathfrak{D}$ or it only contains $\mathfrak{D} \cap X_{p_{0}}.$  If the latter was true, $\Gamma_{g-2}(K_{C} \otimes M_{K_{C}(-p_{0})})$ would be supported on $X_{p_{0}},$ and therefore would have class proportional to $x.$  But this contradicts Lemma \ref{minuspt}.

The same Lemma, together with (\ref{impclass}), shows that $\mathfrak{D}$ and $X_{p_{0}}$ each occur with multiplicity 1 in $\Gamma_{g-2}(K_{C} \otimes M_{K_{C}(-p_{0})})$, thereby establishing (\ref{cyceq}).   
\hfill \qedsymbol

\medskip

Recall that the movable cone of a projective variety $X$ is the closure of the convex cone in $N^{1}_{\mathbb{R}}(X)$ spanned by classes of divisors whose stable base locus has no divisorial component.  Theorem C of \cite{Mus} implies that for a general curve of genus $g \geq 3,$ the class $\theta-x$ of the subordinate locus $\Gamma_{g-1}(K_{C}(-p_{0}))$ lies in the interior of the movable cone of $C_{g-1}$ and has stable base locus $C^{1}_{g-1}.$  In light of this result, the present characterization of $\mathfrak{D}+X_{p_{0}}$ as a higher-rank subordinate locus naturally associated to $K_{C}(-p_{0})$ seems to suggest that the class $(g-2)\theta-(g-1)x$ lies \textit{outside} of the movable cone of $C_{g-2}$ and has stable base locus $\mathfrak{D},$ but we currently have no proof of this.  

\section{The Case Of A Smooth Plane Quintic}
This final section contains the proofs of Theorems VI and VII.  From this point forward, $C$ always denotes a smooth plane quintic, and $\mathcal{Q}$ denotes $\mathcal{O}_{C}(1).$  Note that $\mathcal{Q}^{\otimes 2} \cong K_{C}$ by the adjunction formula. 

Given that the general curve of genus 6 cannot be realized as a smooth plane quintic, it is worth reminding the reader that the classes $x$ and $\theta$ exist in $N^{1}_{\mathbb{R}}(C_{d})$ and are linearly independent for all curves $C$ and all $d \geq 2.$  It can be shown using Proposition 3.9 in \cite{Pir} that $N^{1}_{\mathbb{R}}(C_{d})$ is generated by $x$ and $\theta$ when $C$ is a \textit{general} smooth plane quintic, but we will not need this fact. 

Our first task in this section is to show that $\mathfrak{D}$ is an irreducible effective divisor on $C_{4}.$  The argument for the case of a general curve (which is spread out over the proofs of Proposition III and Theorem I) fails in the plane quintic case for two reasons.  One is that we cannot use the Brill-Noether Theorem to show that $\mathfrak{D}$ is an effective divisor.  Another is that $G^{1}_{5}(C)$ is no longer irreducible; this is a consequence of the next result.  
 
\begin{lem}
\label{bpt}
Let $(L,V)$ be a $\mathfrak{g}^{1}_{5}$ on $C.$
\begin{itemize}
\item[(i)]{If $(L,V)$ is incomplete, then it is a sub-pencil of $|\mathcal{Q}|.$  In particular, the general incomplete $\mathfrak{g}^{1}_{5}$ on $C$ is basepoint-free.}
\item[(ii)]{If $(L,V)$ is complete, then it is of the form $|\mathcal{Q}(r-s)|,$ where $r,s$ are distinct points in $C.$  In particular, every complete $\mathfrak{g}^{1}_{5}$ on $C$ has a base point.}
\end{itemize}
\end{lem}

\begin{proof}
(i) Let $(L,V)$ be an incomplete $\mathfrak{g}^{1}_{5}$ on $C.$  By Clifford's Theorem, the dimension of $|L|$ is exactly 2, and since there is exactly one $\mathfrak{g}^{2}_{5}$ on $C,$ we must have that $L \cong \mathcal{Q}.$ (ii) This follows from  Exercise F-1 in Chapter VIII of \cite{ACGH}. 
\end{proof}

The following characterization of $C^{1}_{4}$ will also be helpful.

\begin{lem}
\label{qsub}
We have the set-theoretic equality $C^{1}_{4}=\Gamma_{4}(\mathcal{Q}).$ 
\end{lem}

\begin{proof}
Any complete linear series of degree 4 on $C$ is necessarily a pencil by Clifford's Theorem.  By Exercise A-12 in Chapter V of \textit{loc. cit.} we have that every $\mathfrak{g}^{1}_{4}$ on $C$ is of the form $|\mathcal{Q}(-s)|$ for some $s \in C;$ the desired equality follows at once.  
\end{proof} 

\begin{prop}
\label{qirr}
$\mathfrak{D}$ is an irreducible effective divisor on $C_{4}.$
\end{prop}
\begin{proof}
We begin by showing that $\mathfrak{D}$ is an effective divisor on $C_{4}$.  The argument from the proof of Proposition III will suffice once we know that $C^{1}_{4}$ is a proper Zariski-closed subset of $C_{4}$ and that $C^{1}_{5}$ is irreducible of the expected dimension.  The first statement follows from Lemma \ref{qsub}.  As for the second statement, Lemma \ref{bpt} implies that the Brill-Noether locus $W^{1}_{5}(C) \subseteq \textnormal{Pic}^{5}(C)$ is a translate of the difference surface $C-C;$ therefore $W^{1}_{5}(C)$ is irreducible of dimension 2.  Since $C^{1}_{5}=a_{5}^{-1}(W^{1}_{5}(C))$ and $a_{5}^{-1}(L) \cong \mathbb{P}^{1}$ for general $L \in W^{1}_{5}(C),$ we have that $C^{1}_{5}$ is irreducible of dimension 3.     

We now proceed to show that $\mathfrak{D}$ is irreducible.  Define the rational map $$\widetilde{\phi}: C \times C \times C \dashrightarrow C_{4}$$ by $(q,r,s) \mapsto q+|\mathcal{Q}(r-s-q)|.$  The degree-4 linear series $|\mathcal{Q}(r-s-q)|$ has nonnegative dimension for all $r,s,q \in C$ since $\mathcal{Q}$ is very ample, and it is zero-dimensional precisely when $r \neq s$ and $r \neq q.$  Therefore $\widetilde{\phi}$ is well-defined on a nonempty Zariski-open subset of $C \times C \times C.$  It remains to show that $\widetilde{\phi}$ dominates $\mathfrak{D}.$

By part (i) of Lemma \ref{bpt} we know that if $D \in C_{4}$ is subordinate to an incomplete $\mathfrak{g}^{1}_{5}$ on $C,$ it is necessarily subordinate to $|\mathcal{Q}|.$  Since $\Gamma_{4}(\mathcal{Q})$ is of codimension 2 in $C_{4}$, it follows that the general element of $\mathfrak{D}$ is subordinate to a complete $\mathfrak{g}^{1}_{5}.$  Lemma \ref{qsub} and part (ii) of Lemma \ref{bpt} imply that if $D \in C_{4}-C^{1}_{4}$ is subordinate to a complete $\mathfrak{g}^{1}_{5}$ on $C,$ then $D+q \in |\mathcal{Q}(r-s)|$ for some $r,s,q \in C$ such that $r \neq s$ and $r \neq q.$  In particular, $D=\widetilde{\phi}(q,r,s).$
\end{proof}

For the proof of Theorem VI, we will need the following Lemma, whose derivation is analogous to that of Lemma \ref{minuspt}.

\begin{lem}
\label{qvir}
The locus $\Gamma_{4}(K_{C} \otimes M_{\mathcal{Q}})$ is a virtual divisor on $C_{4}$ whose class is $2\theta-3x.$ \hfill \qedsymbol
\end{lem}

\subsection{\textnormal{\textit{Proof of Theorem VI}}}
By (\ref{impclass}) and Lemma \ref{qvir}, the equality of cycles (\ref{quin}) will be established once we check equality of the underlying sets.  Since $\Gamma_{4}(K_{C} \otimes M_{\mathcal{Q}})$ is a virtual divisor and $\mathfrak{D}$ is an irreducible effective divisor by Proposition \ref{qirr}, it suffices to show that $\Gamma_{4}(K_{C} \otimes M_{\mathcal{Q}})$ is a nonempty subset of $\mathfrak{D}.$  This will be accomplished by a useful characterization of $\Gamma_{4}(K_{C} \otimes M_{\mathcal{Q}}),$ which we now describe.

Let $\iota_{\mathcal{Q}}:C \hookrightarrow \mathbb{P}^{2}$ be the embedding induced by $\mathcal{Q}.$  Recall that for a closed subscheme $X$ of $\mathbb{P}^{N},$ the space $H^{0}(\mathbb{P}^{N},\Omega^{1}_{\mathbb{P}^{N}}(1) \otimes \mathcal{I}_{X}(k))$ is the kernel of the multiplication map $$H^{0}(\mathbb{P}^{N},\mathcal{O}_{\mathbb{P}^{N}}(1)) \otimes H^{0}(\mathbb{P}^{N},\mathcal{I}_{X}(k)) \rightarrow H^{0}(\mathbb{P}^{N},\mathcal{I}_{X}(k+1))$$ i.e. the space of linear syzygies among the degree-$k$ elements of the saturated homogenous ideal of $X.$  As mentioned in the Introduction, for each $D \in C_{4}$ the space $H^{0}(C,K_{C} \otimes M_{\mathcal{Q}}(-D))$ is the kernel of the multiplication map $$H^{0}(C,\mathcal{Q}) \otimes H^{0}(C,K_{C}(-D)) \rightarrow H^{0}(C,\mathcal{Q} \otimes K_{C}(-D)).$$  Restriction from $\mathbb{P}^{2}$ to $C$ yields the following commutative diagram with exact rows. 
\vskip20pt
\begin{tiny}
$\begin{CD}
0  @>>> {H^{0}(\mathbb{P}^{2},\Omega^{1}_{\mathbb{P}^{2}}(1) \otimes \mathcal{I}_{\iota_{\mathcal{Q}}(D)}(2))}              @>>> {H^{0}(\mathbb{P}^{2},\mathcal{O}_{\mathbb{P}^{2}}(1)) \otimes H^{0}(\mathbb{P}^{2},\mathcal{I}_{\iota_{\mathcal{Q}}(D)}(2))}    @>>> {H^{0}(\mathbb{P}^{2},\mathcal{I}_{\iota_{\mathcal{Q}}(D)}(3))}\\ 
@.     @VVV                                                @VVV																	 @VVV\\
0  @>>> H^{0}(C,E^{\ast}_{\mathcal{Q}} \otimes K_{C}(-D))     @>>> {H^{0}(C,\mathcal{Q}) \otimes H^{0}(K_{C}(-D))} @>>> {H^{0}(C,K_{C} \otimes \mathcal{Q}(-D))}\\ 
\end{CD}$
\end{tiny}
\vskip20pt
Since $h^{i}(\mathbb{P}^{2},\mathcal{O}_{\mathbb{P}^{2}}(-3))=h^{i}(\mathbb{P}^{2},\mathcal{O}_{\mathbb{P}^{2}}(-2))=0$ for $i=0,1,$ the right two vertical arrows are isomorphisms, and therefore so is the leftmost vertical arrow.  It follows that $D \in \Gamma_{4}(K_{C} \otimes M_{\mathcal{Q}})$ precisely when the quadrics in $\mathbb{P}^{2}$ which vanish along the subscheme $\iota_{\mathcal{Q}}(D)$ admit linear syzygies.

For any $D \in C_{4},$ there are three distinct possibilities for $\iota_{\mathcal{Q}}(D):$ 
\begin{itemize}
\item[(a)]{$\iota_{\mathcal{Q}}(D)$ is collinear.}
\item[(b)]{$\iota_{\mathcal{Q}}(D)$ is noncollinear, but contains a collinear subscheme of length 3.}
\item[(c)]{$\iota_{\mathcal{Q}}(D)$ contains no collinear subscheme of length 3.}
\end{itemize}  
It follows from the discussion in Example 3.12 of \cite{Eis} and the upper-semicontinuity of Betti numbers that the quadrics in $\mathbb{P}^{2}$ which vanish along $\iota_{\mathcal{Q}}(D)$ admit linear syzygies precisely when (a) or (b) holds.  In particular, $\Gamma_{4}(\mathcal{Q})=C^{1}_{4}$ is a subset of $\Gamma_{4}(K_{C} \otimes M_{\mathcal{Q}}),$ so that $\Gamma_{4}(K_{C} \otimes M_{\mathcal{Q}})$ is nonempty.

It remains to show that $\Gamma_{4}(K_{C} \otimes M_{\mathcal{Q}}) \subseteq \mathfrak{D}.$  If (a) holds, and $\ell \subseteq \mathbb{P}^{2}$ is the line spanned by $\iota_{\mathcal{Q}}(D),$  then $D$ is subordinate to the (incomplete) $\mathfrak{g}^{1}_{5}$ on $C$ obtained by projecting away from a point on $\ell$ which does not lie on $C.$  If (b) holds, then $D=D'+r,$ where $\iota_{\mathcal{Q}}(D')$ is collinear of length 3 and $r$ is not on the line spanned by $\iota_{\mathcal{Q}}(D').$  If $s \in C$ is one of the points for which $D' \in \Gamma_{3}(\mathcal{Q}(-s)),$ we have that $D \in \Gamma_{4}(\mathcal{Q}(r-s)).$  Since $|\mathcal{Q}(r-s)|$ is a $\mathfrak{g}^{1}_{5},$ we are done.    \hfill \qedsymbol

\subsection{\textnormal{\textit{Proof of Theorem VII}}}
We begin by determining the components of the subordinate locus associated to a complete $\mathfrak{g}^{1}_{5}$ on $C.$

\begin{prop}
\label{comps}
For $p,s \in C, p \neq s,$ define $Z_{p,s}:=p+\Gamma_{3}(\mathcal{Q}(-s)).$  Then we have the equality of cycles  
\begin{equation}
\label{pieces}
\Gamma_{4}(\mathcal{Q}(p-s))=Z_{p,s}+\Gamma_{4}(\mathcal{Q}(-s))
\end{equation}
\end{prop} 

\begin{proof}
We first establish the equality at the level of sets and then use a specialization argument to show that each component occurs with multiplicity 1.  Let $D \in \Gamma_{4}(\mathcal{Q}(p-s))$ be given.  Then there exists $r' \in C$ such that $D+s+r' \in |\mathcal{Q}(p)|,$ and since $|\mathcal{Q}(p)|$ is a $\mathfrak{g}^{2}_{6}$ with base point $p$ by Clifford's Theorem, we have that either $r'=p$ or $D-p$ is effective.  In the first case, we have that $D \in \Gamma_{4}(\mathcal{Q}(-s))$ and in the second case we have that $D \in Z_{p,s}.$  The reverse inclusion $Z_{p,s} \cup \Gamma_{4}(\mathcal{Q}(-s)) \subseteq \Gamma_{4}(\mathcal{Q}(p-s))$ follows immediately from the definitions.

By Lemma \ref{genus}, $Z_{p,s} \cong C.$  Furthermore, $\Gamma_{4}(\mathcal{Q}(-s)) \cong \mathbb{P}^{1},$ and since $Z_{p,s}$ is a smooth subvariety of the divisor $X_{p}$ which meets $\Gamma_{4}(\mathcal{Q}(-s))$ transversally in one point by Corollary \ref{exc}, we have that $Z_{p,s}$ meets $\Gamma_{4}(\mathcal{Q}(-s))$ properly in one point.  Therefore the reducible curve $Z_{p,s} + \Gamma_{4}(\mathcal{Q}(-s))$ has arithmetic genus 6.  Another application of Lemma \ref{genus} shows that $\Gamma_{4}(\mathcal{Q},V) \cong C$ whenever $(\mathcal{Q},V)$ is a basepoint-free pencil of degree 5.  Since $G^{1}_{5}(C)$ is connected (this follows, for instance, from Theorem \ref{fullaz}) we have that the arithmetic genus of $\Gamma_{4}(\mathcal{Q}(p-s))$ is also 6, so that $(\ref{pieces})$ holds.    
\end{proof}

The next two propositions conclude the proof of Theorem VII.  We denote by $\mathfrak{Z}$ the numerical class of the curves $Z_{p,s}.$

\begin{prop}
\label{empty}
$(\theta-2x) \cdot \mathfrak{Z}=0.$
\end{prop}

\begin{proof}
Combining Proposition \ref{comps} and Corollary \ref{exc} with our calculations from the proof of Theorem I, we have that 
\begin{align*}
\theta \cdot \mathfrak{Z} &= \theta \cdot \gamma_{4}(\mathfrak{g}^{1}_{5})-\theta \cdot \gamma_{4}(\mathfrak{g}^{1}_{4})=6\\
x \cdot \mathfrak{Z} &= x \cdot \gamma_{4}(\mathfrak{g}^{1}_{5})-x \cdot \gamma_{4}(\mathfrak{g}^{1}_{4})=3
\end{align*}
and the result follows.
\end{proof}

\begin{prop}
If $\mathfrak{E}$ is an effective divisor on $C_{4}$ whose class is proportional to $\theta-tx$ for some $t>0,$ then $\mathfrak{E} \cdot \mathfrak{Z} \geq 1.$
\end{prop}
\begin{proof}
We may suppose without loss of generality that $\mathfrak{E}$ is irreducible.  The conclusion follows immediately from the calculations in the previous proof if $t=\frac{3}{2},$ so we will assume otherwise from now on.  In particular, we are assuming that $\mathfrak{E} \neq \Gamma_{4}(K_{C} \otimes E^{\ast}_{\mathcal{Q}}).$  Since the curves $Z_{p,s}$ cover $\Gamma_{4}(K_{C} \otimes E^{\ast}_{\mathcal{Q}}),$ there necessarily exist distinct $p',s' \in C$ such that $Z_{p',s'} \nsubseteq \mathfrak{E}.$  Therefore we have that $\mathfrak{E} \cdot \mathfrak{Z} \geq 0,$ and it remains to check that this intersection is positive.

By Nakamaye's theorem on stable base loci of nef and big divisors (e.g. Theorem 10.3.5 in \cite{Laz}) applied to the class $\theta,$ the numerical hypothesis on $\mathfrak{E}$ implies that $C^{1}_{4} \subset \mathfrak{E}.$ Since $Z_{p',s'} \cap C^{1}_{4} \neq \emptyset,$ the result follows. 
\end{proof}

\vskip15pt

\textit{Funding:} This work was partially supported by the National Science Foundation [RTG DMS-0502170]. 

\vskip15pt

\textit{Acknowledgments:} The author would like to thank S. Casalaina-Martin, E. Cotterill, J. Kass, R. Lazarsfeld, M. Nakamaye, G. Pacienza, J. Ross, and J. Starr for valuable discussions and correspondence.  He would also like to thank the anonymous referees for numerous helpful suggestions and for pointing towards a gap in the initial statement of Proposition III.   

\medskip

\noindent{Department of Mathematics, University of Michigan\\
East Hall, 530 Church Street\\
Ann Arbor, MI 48109-1043\\
email: $\tt{ymustopa@umich.edu}$}


\begin{thebibliography}{10}
\vspace{0.1in}
\bibitem[ACGH]{ACGH}
\begin{singlespace}
Arbarello, Enrico, Maurizio Cornalba, Phillip Griffiths and Joe Harris.
\newblock {\sl Geometry of Algebraic Curves, Volume I}, 
\newblock New York: Springer-Verlag, 1985.
\end{singlespace}
\begin{singlespace}
\bibitem[ApFa]{ApFa}
Aprodu, Marian, and Gavril Farkas.
\newblock {``Koszul cohomology and applications to moduli"},
\newblock to appear in {\sl Aspects of vector bundles and moduli, Clay Mathematical Institute Volume 10}
\end{singlespace}
\begin{singlespace}
\bibitem[Chan]{Chan}
Chan, Kungho.
\newblock {``A characterization of double covers of curves in terms of the ample cone of second symmetric product"},
\newblock {\sl J. Pure. App. Alg.} 212 no. 12 (2008): 2623-32. 
\end{singlespace}
\begin{singlespace}
\bibitem[Cot]{Cot}
Cotterill, Ethan.  
\newblock {``Geometry of curves with exceptional secant planes: linear series along the general curve"},
\newblock to appear in {\sl Math. Zeit.}
\end{singlespace}
\begin{singlespace}
\bibitem[Deb]{Deb}
Debarre, Olivier.
\newblock {``Seshadri Constants Of Abelian Varieties"},
\newblock in {\sl The Fano Conference}, Univ. Torino, Turin, 2004: 379-94.
\end{singlespace}
\begin{singlespace}
\bibitem[Ein]{Ein}
Ein, Lawrence.
\newblock {``A remark on the syzygies of the generic canonical curves"},
\newblock {\sl J. Differential Geom.} 26 (1987): 361-65. 
\end{singlespace}
\begin{singlespace}
\bibitem[Eis]{Eis}
Eisenbud, David.
\newblock {\sl The Geometry of Syzygies},
\newblock New York: Springer-Verlag, 2005.
\end{singlespace}
\begin{singlespace}
\bibitem[Fa]{Fa}
Farkas, Gavril.
\newblock {``Higher ramification and varieties of secant divisors on the generic curve"},
\newblock {\sl J. London Math. Soc.} 78 (2008): 418-40.
\end{singlespace}
\begin{singlespace}
\bibitem[vdGK]{vdGK}
van der Geer, Gerard, and Alexis Kouvidakis.
\newblock {``Cycle Relations On Jacobian Varieties"},
\newblock {\sl Compositio Math} 143 (2007): 900-08.
\end{singlespace}
\begin{singlespace}
\bibitem[Gr]{Gr}
Green, Mark.
\newblock {``Koszul Cohomology And The Geometry Of Projective Varieties"},
\newblock {\sl J. Differential Geom.} 19 (1984): 125-71. 
\end{singlespace}
\begin{singlespace}
\bibitem[Herb]{Herb}
Herbaut, Fabien.
\newblock {``Algebraic Cycles On The Jacobian Of A Curve With A Linear System Of Given Dimension"},
\newblock {\sl Compositio Math.} 143 (2007): 883-99.
\end{singlespace}
\begin{singlespace}
\bibitem[Iz1]{Iz1}
Izadi, Elham.
\newblock {``Deforming Curves In Jacobians To Non-Jacobians I: Curves in $C^{(2)}$"},
\newblock {\sl Geom. Dedicata} 116 (2005): 87-109.
\end{singlespace}
\begin{singlespace}
\bibitem[Iz2]{Iz2}
Izadi, Elham.
\newblock {``Deforming curves In Jacobians To Non-Jacobians II. Curves in $C^{(e)}, 3 \leq e \leq g-3$"},
\newblock {\sl Geom. Dedicata} 115 (2005): 33-63.
\end{singlespace}
\begin{singlespace}
\bibitem[Kou]{Kou}
Kouvidakis, Alexis.
\newblock {``Divisors on Symmetric Products of Curves"},
\newblock {\sl Trans. Amer. Math. Soc.} 337 (1993): 117-28.
\end{singlespace}
\begin{singlespace}
\bibitem[Laz]{Laz}
Lazarsfeld, Robert.
\newblock \sl{Positivity in Algebraic Geometry II},
\newblock New York: Springer-Verlag, 2004.
\end{singlespace}
\begin{singlespace}
\bibitem[Mus]{Mus}
Mustopa, Yusuf.
\newblock {``Residuation of Linear Series and the Effective Cone of $C_{d}$"},
\newblock to appear in {\sl Amer. J. Math.}
\end{singlespace}
\begin{singlespace}
\bibitem[Pac]{Pac}
Pacienza, Gianluca.
\newblock {``On the nef cone of symmetric products of a generic curve"},
\newblock {\sl Amer. J. Math.} 125 (2003): 1117-35.
\end{singlespace}
\begin{singlespace}
\bibitem[PR]{PR}
Paranjape, Kapil, and Sundraraman Ramanan.
\newblock {``On the canonical ring of a curve"},
\newblock in {\sl Algebraic Geometry and Commutative Algebra II}, Kinokuniya, Tokyo, 1984. 
\end{singlespace}
\begin{singlespace}
\bibitem[Pir]{Pir}
Pirola, Gian Pietro.
\newblock {``Base Number Theorem for Abelian Varieties"},
\newblock {\sl Math. Ann.} 282 (1988): 361-368
\end{singlespace}
\begin{singlespace}
\bibitem[Tei]{Tei}
Teixidor i Bigas, Montserrat.
\newblock {``Syzygies using vector bundles"},
\newblock {\sl Trans. Amer. Math. Soc.} 359 (2007): 897-908.
\end{singlespace}
\begin{singlespace}
\bibitem[V]{V}
Voisin, Claire.
\newblock {``Green's generic syzygy conjecture for curves of even genus lying on a K3 surface"},
\newblock {\sl J. Eur. Math. Soc.} 2 (2004): 363-404.
\end{singlespace}
\end{thebibliography}
\end{document}